\documentclass{article}
\usepackage{amssymb}
\usepackage{amsmath}

\newtheorem{theorem}{Theorem}

\newtheorem{definition}{Definition}

\newtheorem{proposition}{Proposition}

\newenvironment{proof}[1][Proof]{\noindent\textbf{#1.} }{\ \rule{0.5em}{0.5em}}
\input{tcilatex}

\begin{document}

\title{The problem of the classification of the nilpotent class $2$ torsion
free groups up to the geometrically equivalence.}
\author{{\Large A. Tsurkov} \\
%EndAName
\textit{Jerusalem College of Technology,}\\
\textit{P. B. 36177, Jerusalem, 91361, Israel.}\\
tsurkov@jct.ac.il\bigskip \\
\textit{The Open University of Israel.}}
\maketitle

\begin{abstract}
In this paper we consider the problem of classification of the nilpotent
class $2$ finitely generated torsion free groups up to the geometric
equivalence. By a very easy technique it is proved that this problem is
equivalent to the problem of classification of the complete (in the Maltsev
sense) nilpotent torsion free finite rank groups up to the isomorphism. This
result, allows us to once more comprehend the complication of the problem of
the classification of the quasi-varieties of nilpotent class $2$ groups. It
is well known that the variety of a nilpotent class $s$ (for every $s\in 
%TCIMACRO{\U{2115} }%
%BeginExpansion
\mathbb{N}
%EndExpansion
$) groups is Noetherian. So the problem of the classification of the
quasi-varieties generated even by a single nilpotent class $2$ finitely
generated torsion free group, is equivalent to the problem of classification
of complete (in the Maltsev sense) nilpotent torsion free finite rank groups
up to the isomorphism.
\end{abstract}

\section{Historical review and methodology.}

The one of the question which naturaly appears in all algebraic studies, is
the question of the classification (up to isomorphism) of the algebraic
objects from some class. One of the classical example of this kind of
results is the classification of the semisimple finite dimensional
associative algebras over fields. Also we have the very easely observed
classification of the finitely generated abelian groups. Both of these
classification were achieved many years ago. About 60 years ago, the
classification of the simple Lie algebras over $%
%TCIMACRO{\U{2102} }%
%BeginExpansion
\mathbb{C}
%EndExpansion
$ was achieved. The classification of the finite simple groups, is a newer
result of this kind. This result requested a huge effort of many
mathematicians.

The natural next step after the classification of the finitely generated
abelian groups, is the classification of the finitely generated nilpotent
class $2$ groups, in particular the classification of the finitely generated
torsion free nilpotent class $2$ groups. In the 1970s the research in this
area was very active. We can remember \cite{GrSch}, \cite{GrSeSt}, auxiliary
technical work \cite{Sch} and others. Achievements were summarized in \cite%
{GrSe}. A full classification up to isomorphism was established only for
finitely generated torsion free nilpotent class $2$ groups of Hirsh length $%
6 $. Such modest results indicate the complication of the problem. This
complication follows, as it will be seen from our survey, first of all, from
the complication of the wild matrix problem.

Also in the 1970s, the problem of classification of the nilpotent class $2$ $%
p$-groups up to isomorphism was considered in \cite{Serg1}. It was proved
that this problem can be reduced to the wild matrix problem even when the
rank of the center of the groups is equal to $2$.

It is known that for every nilpotent torsion free group $G$ there is a
Maltsev completion $\sqrt{G}$ - the minimal group, such that $G\subset \sqrt{%
G}$ and for every $x\in \sqrt{G}$ and every $n\in \mathbf{%
%TCIMACRO{\U{2115} }%
%BeginExpansion
\mathbb{N}
%EndExpansion
}$ there exists $x^{\frac{1}{n}}\in \sqrt{G}$, such that $\left( x^{\frac{1}{%
n}}\right) ^{n}=x$. The element $x^{\frac{1}{n}}\in \sqrt{G}$ is uniquely
defined by $x\in \sqrt{G}$ and $n\in \mathbf{%
%TCIMACRO{\U{2115} }%
%BeginExpansion
\mathbb{N}
%EndExpansion
}$. The $\sqrt{G}$ is the nilpotent group of the same class of nilpotence as 
$G$. It is clear that if two nilpotent torsion free groups are isomorphic,
then their Maltsev completions are isomorphic too. The inverse, of course,
is false. So the classification of the complete (in the Maltsev sense)
nilpotent torsion free groups of the finite rank up to isomorphism is a
simpler problem than the classification up to isomorphism of the arbitrary
finitely generated nilpotent class $2$ groups. But even in this simpler
problem we have (see \cite{GrSe}) a solution, only in the case when the rank
of the center of the groups is not greater than $2$.

If the problem of the classification up to isomorphism is so complicated, it
is natural to consider a less delicate classification. The notion of the
geometric equivalence of universal algebras (see \cite{Pl1}), was
investigated in 1995 by B. Plotkin. By \cite{Pl2} two finitely generated
universal algebras $A_{1}$ and $A_{2}$ from some variety $\Theta $ are
geometrically equivalent if and only if, the first one can be embedded into
some direct power of the second and vice versa (we denote\ $A_{1}\sim A_{2}$%
). So, the classification up to geometric equivalence is less delicate then
the classification up to isomorphism.

Classification of the nilpotent groups up to geometric equivalence is
especialy interesting because in the case of the nilpotent groups geometric
equivalence is closely connected with the logic proprieties of groups: two
nilpotent groups\ are geometrically equivalent if and only if, they have the
same quasi-identities (see \cite{Ts}). So classification of nilpotent
groups, up to geometric equivalence is an equivalent to the classification
of the quasi-varieties generated by a single nilpotent group.

The classification of the finitely generated abelian groups up to geometric
equivalence was achieved in \cite{Be}. It was proved that two abelian groups
are geometrically equivalent if and only if, for every prime number $p$ the
exponents of their corresponding $p$-Sylow subgroups coincide, and if one of
these group is not periodic, then the second group is not periodic either.
So the classification of finitely generated abelian groups up to geometric
equivalence, is in principal, simpler than the classification of these
groups up to isomorphism. The classification of the torsion free abelian
groups up to geometric equivalence is trivial: all these groups are
geometrically equivalent.

In the case of the nilpotent class $2$ groups we have a different situation.
Even classification of the finitely generated torsion free nilpotent class $%
2 $ groups up to geometric equivalence is a very complicated problem.

By \cite[Theorem 1]{Ts} every finitely generated torsion free nilpotent
class $2$ group geometrically equivalent to it's Maltsev completion. So for
resolving of our problem it is enough to classify up to geometric
equivalence the nilpotent class $2$ finite rank torsion free complete groups.

It is well known (see for example \cite[Chapter 8]{Ba}) that in every
nilpotent Lie $\mathbf{%
%TCIMACRO{\U{211a} }%
%BeginExpansion
\mathbb{Q}
%EndExpansion
}$ -algebra $L$, we can define multiplication by Campbell-Hausdorff formula.
With this multiplication $L$ will be a group, which we denote $L^{\circ }$.
The group $L^{\circ }$ will be torsion free and complete. It has the same
class of nilpotency as algebra $L$. Conversely for every complete nilpotent
torsion free group $A$, there is nilpotent Lie $\mathbf{%
%TCIMACRO{\U{211a} }%
%BeginExpansion
\mathbb{Q}
%EndExpansion
}$-algebra $L$, such that $A\cong L^{\circ }$. Algebra $L$ has the some
class of nilpotency as group $A$. The homomorphisms (epimorphisms,
monomorphisms, isomorphisms) of the nilpotent Lie $\mathbf{%
%TCIMACRO{\U{211a} }%
%BeginExpansion
\mathbb{Q}
%EndExpansion
}$-algebras coincide with the homomorphisms (epimorphisms, monomorphisms,
isomorphisms) of the corresponding groups and vice versa. In other words,
the functor $\Gamma :L\rightarrow L^{\circ },\left( \lambda
:L_{1}\rightarrow L_{2}\right) \rightarrow \left( \lambda :L_{1}^{\circ
}\rightarrow L_{2}^{\circ }\right) $ provides an isomorphism from the
category of the nilpotent class $s$ Lie $\mathbf{%
%TCIMACRO{\U{211a} }%
%BeginExpansion
\mathbb{Q}
%EndExpansion
}$-algebras to the category of the nilpotent class $s$ torsion free complete
groups. By this isomorphism to the finite dimensional Lie $\mathbf{%
%TCIMACRO{\U{211a} }%
%BeginExpansion
\mathbb{Q}
%EndExpansion
}$ -algebra corresponds the nilpotent group of the finite rank and vice
versa.

So two complete nilpotent torsion free finite rank groups $%
A_{1}=L_{1}^{\circ }$ and $A_{2}=L_{2}^{0}$, are isomorphic if and only if,
the Lie $\mathbf{%
%TCIMACRO{\U{211a} }%
%BeginExpansion
\mathbb{Q}
%EndExpansion
}$-algebras $L_{1}$ and $L_{2}$ are isomorphic. And two complete nilpotent
torsion free finite rank groups $A_{1}=L_{1}^{\circ }$ and $%
A_{2}=L_{2}^{\circ }$ are geometrically equivalent if and only if, the Lie $%
\mathbf{%
%TCIMACRO{\U{211a} }%
%BeginExpansion
\mathbb{Q}
%EndExpansion
}$-algebras $L_{1}$ and $L_{2}$ are geometrically equivalent, i.e. as it was
stated above, if and only if the algebra $L_{1}$ can be embedded into some
direct power of the algebra $L_{2}$ and vice versa. So for researching our
problem, we can concentrate on the geometric equivalence of the finite
dimension nilpotent class $2$ Lie $\mathbf{%
%TCIMACRO{\U{211a} }%
%BeginExpansion
\mathbb{Q}
%EndExpansion
}$-algebras. It will be proved in this paper, that the problem of
classification of the finite dimension nilpotent class $2$ Lie $\mathbf{%
%TCIMACRO{\U{211a} }%
%BeginExpansion
\mathbb{Q}
%EndExpansion
}$-algebras up to the geometric equivalence, is equivalent to the problem of
classification of these algebras up to the isomorphism. It means that the
problem of the classification of the nilpotent class $2$ finite rank torsion
free complete groups up to the geometric equivalence, is equivalent to the
problem of classification of these groups up to the isomorphism.

The problem of classification of the finite dimension nilpotent class $2$
Lie algebras over an algebraic closed field up to the isomorphism was
considered in \cite{Serg2} and \cite{BLS}. In \cite{Serg2} this problem was
resolved when the dimension of the center of the algebra is not great then $%
2 $. In \cite{BLS} it was proved that the problem of classification of the
finite dimension nilpotent class $2$ Lie algebras over an algebraic closed
field up to the isomorphism when the dimension of the center of the algebra
is great then $2$ is equivalent to the wild problem.

In \cite{Be} it was proved that if $A_{1}\sim A_{2}$ and $B_{1}\sim B_{2}$ ($%
A_{1},A_{2},B_{1},B_{2}$ are arbitrary universal algebras from some variety $%
\Theta $) then $A_{1}\oplus B_{1}\sim A_{2}\oplus B_{2}$. So for
classification up to the geometric equivalence it is enough to consider the
algebras which are can not be decomposed to the direct sum. In our
situation, we can consider only Lie $\mathbf{%
%TCIMACRO{\U{211a} }%
%BeginExpansion
\mathbb{Q}
%EndExpansion
}$-algebras $L$ which fulfill 
\begin{equation}
\left[ L,L\right] =Z\left( L\right) ,  \label{com_cond}
\end{equation}%
where $Z\left( L\right) $ is a center of the algebra $L$.

The Lie $\mathbf{%
%TCIMACRO{\U{211a} }%
%BeginExpansion
\mathbb{Q}
%EndExpansion
}$-algebra $L$, which fulfills this condition, can be considered as the
direct sum of the $\mathbf{%
%TCIMACRO{\U{211a} }%
%BeginExpansion
\mathbb{Q}
%EndExpansion
}$-linear spaces $L=V\oplus W$ (form this place and below we considered the
direct sum only in the category of the $\mathbf{%
%TCIMACRO{\U{211a} }%
%BeginExpansion
\mathbb{Q}
%EndExpansion
}$-linear spaces), where $W=Z\left( L\right) $, $V\cong L/Z\left( L\right) $%
. The Lie brackets in $L$, define the skew symmetric non singular bilinear
mapping $\omega _{L}:V\times V\ni \left( v_{1},v_{2}\right) \rightarrow %
\left[ v_{1},v_{2}\right] \in W$. (For arbitrary skew symmetric bilinear
mapping $\omega :V\times V\rightarrow W$ we denote $\ker \omega =\left\{
x\in V\mid \forall v\in V\left( \omega \left( x,v\right) =0\right) \right\} $
and we say that $\omega $ is singular if $\ker \omega \neq \left\{ 0\right\} 
$; other skew symmetric bilinear mappings we call non singular.)

Contrariwise, if we have two $\mathbf{%
%TCIMACRO{\U{211a} }%
%BeginExpansion
\mathbb{Q}
%EndExpansion
}$-linear spaces $V$ and $W$ and the skew symmetric bilinear mapping $\omega
:V\times V\rightarrow W$, then in the direct sum $L=V\oplus W$ we can define
the Lie brackets: $\left[ v_{1}+w_{1},v_{2}+w_{2}\right] =\omega \left(
v_{1},v_{2}\right) $ ($v_{1},v_{2}\in V$, $w_{1},w_{2}\in W$). If $\omega $
is a non singular, then $Z\left( L\right) =W$ and condition (\ref{com_cond})
is an equivalent to the condition%
\begin{equation}
\omega \left( V,V\right) =W.  \label{im_cond}
\end{equation}%
If $L=V\oplus W$ and $\dim V=n$, $\dim W=m$, $\left\{ v_{1},\ldots
,v_{n}\right\} $ is a basis of $V$, $\left\{ w_{1},\ldots ,w_{m}\right\} $
is a basis of $W$, then the skew symmetric bilinear mapping $\omega $ is
defined by $m$ skew symmetric matrices of the size $n\times n$: $A^{\left(
1\right) },\ldots ,A^{\left( m\right) }$, such that $\left[ v_{i},v_{j}%
\right] =\sum\limits_{k=1}^{m}a_{i,j}^{\left( k\right) }w_{k}$ ($1\leq
i,j\leq n$). There is a homomorphism of the Lie algebras with condition (\ref%
{com_cond}) $\lambda :L=V_{L}\oplus W_{L}\rightarrow S=V_{S}\oplus W_{S}$ if
and only if, there is a pair of the linear mappings $\left( \varphi ,\psi
\right) $ such that $\varphi :V_{L}\rightarrow V_{S}$, $\psi
:W_{L}\rightarrow W_{S}$ and $\omega _{S}\left( \varphi \left( v_{1}\right)
,\varphi \left( v_{2}\right) \right) =\psi \omega _{L}\left(
v_{1},v_{2}\right) $ for every $v_{1},v_{2}\in V_{L}$. $\lambda =\varphi
\oplus \psi $ holds.

By using this approach in \cite{Ts} (and by using the isomorphism from the
category of the nilpotent class $2$ Lie $\mathbf{%
%TCIMACRO{\U{211a} }%
%BeginExpansion
\mathbb{Q}
%EndExpansion
}$ -algebras to the category of the nilpotent class $2$ torsion free
complete groups) the problem of the classification of the nilpotent class $2$
finite rank torsion free groups, whose centers have rank no more then $2$,
up to the geometric equivalence was deeply researched. It was proved two
theorems:

\begin{enumerate}
\item Theorem 3. Two nilpotent torsion free class $2$ finitely generated
groups $G_{1}$ and $G_{2}$ with the cyclic center are geometrically
equivalent ($G_{1}\sim G_{2}$) if and only if their Maltsev completions are
isomorphic: $\sqrt{G_{1}}\cong \sqrt{G_{2}}$.

\item Theorem 4. Let $G_{1},G_{2}$ two nilpotent torsion free class $2$
finitely generated groups, whose centers have rank $2$. Then $G_{1}\sim
G_{2} $ if and only if, or there is a nilpotent torsion free class $2$
finitely generated group with the cyclic center $N$, such that $G_{1}\sim
N\sim G_{2}$, or $\sqrt{G_{1}}\cong \sqrt{G_{2}}$.
\end{enumerate}

Also Proposition 1 and Proposition 2 from \cite{Ts}, which formulated by the
language of the properties of the skew symmetric bilinear forms, provide us
tools for finding when for two nilpotent torsion free class $2$ finitely
generated groups $G_{1}$ and $G_{2}$, whose centers have rank $2$, fulfills
the first or the second condition of the Theorem 4.

\section{New results.}

Bellow the word "algebra" means: nilpotent class $2$ finite dimension Lie $%
\mathbf{%
%TCIMACRO{\U{211a} }%
%BeginExpansion
\mathbb{Q}
%EndExpansion
}$-algebra.

\begin{definition}
\label{decomp}We say that the algebra is \textit{geometrically decomposable}
if it is geometrically equivalent to the direct product of some of its
nontrivial subalgebras. Other algebras we call \textit{geometrically
indecomposable}.
\end{definition}

\begin{proposition}
\label{equivalent_isom}\textit{If two geometrically indecomposable algebras
are geometrically equivalent, then they are isomorphic.}
\end{proposition}

\begin{proof}
We assume that $L$ and $S$ are geometrically indecomposable and $L$
geometrically equivalent to $S$. There is a family of
homomorphisms\linebreak\ $\left\{ \lambda _{i}:L\rightarrow S\mid i\in
I\right\} $ such that $\bigcap\limits_{i\in I}\ker \lambda _{i}=\left\{
0\right\} $. Also exists a family of homomorphisms $\left\{ \sigma
_{j}:S\rightarrow L\mid j\in J\right\} $ such that $\bigcap\limits_{j\in
J}\ker \sigma _{j}=\left\{ 0\right\} $. If $L=\left\{ 0\right\} $ then $%
S=\left\{ 0\right\} $ and vice versa, so we can assume that $L,S\neq \left\{
0\right\} $.

We assume that $\ker \lambda _{i}\neq \left\{ 0\right\} $ for every $i\in I$%
. We consider the family of endomorphisms $\left\{ \sigma _{j}\lambda
_{i}:L\rightarrow L\mid j\in J,i\in I\right\} $. $\bigcap\limits_{\substack{ %
j\in J  \\ i\in I}}\ker \sigma _{j}\lambda _{i}=\left\{ 0\right\} $, so
there exists an embedding $L\hookrightarrow \prod\limits_{\substack{ j\in J 
\\ i\in I}}\mathrm{im}\sigma _{j}\lambda _{i}$. $\ker \sigma _{j}\lambda
_{i}\supset \ker \lambda _{i}\neq \left\{ 0\right\} $, so, by reason of
dimensions, $\mathrm{im}\sigma _{j}\lambda _{i}$ is not equal to $L$. We
have that $L$ is geometrically equivalent to $\prod\limits_{\substack{ j\in
J  \\ i\in I}}\mathrm{im}\sigma _{j}\lambda _{i}$. If $L\neq \left\{
0\right\} $ then there is a nonzero factor in this product. This is a
contradiction.

By symmetry we achieve a contradiction when we assume that $\ker \sigma
_{j}\neq \left\{ 0\right\} $ for every $j\in J$.

So there exist $i\in I$, such that $\ker \lambda _{i}=\left\{ 0\right\} $
and $j\in J$ , such that $\ker \sigma _{j}=\left\{ 0\right\} $. By reason of
dimensions, $L$ and $S$ are isomorphic.
\end{proof}

\begin{proposition}
\label{decompos_cond}\textit{If the algebra }$L=V\oplus W$\textit{\ is
geometrically decomposable then there exists a family of linear mappings }$%
\left\{ \psi _{i}:W\rightarrow W\mid i\in I\right\} $\textit{\ such that all
skew symmetric bilinear mappings }$\psi _{i}\omega $\textit{\ are singular
and }$\bigcap\limits_{i\in I}\ker \psi _{i}=\left\{ 0\right\} $ ($\omega
:V\times V\rightarrow W$ is a skew symmetric bilinear mapping defined by the
Lie brackets of the \textit{algebra }$L$)\textit{.}
\end{proposition}

\begin{proof}
Let the algebra $L$ be geometrically decomposable. Then $L$ is geometrically
equivalent to $\prod\limits_{i\in I_{0}}L_{i}$, where $L_{i}$ ($i\in I_{0}$)
is nontrivial subalgebras of $L$. So there exists an embedding $%
L\hookrightarrow \left( \prod\limits_{i\in I_{0}}L_{i}\right) ^{J}$. But we
can write $\left( \prod\limits_{i\in I_{0}}L_{i}\right)
^{J}=\prod\limits_{i\in I}L_{i}$,where $I=I_{0}\times J$ and $L_{i}$ ($i\in
I $) is also nontrivial subalgebras of $L$. By Remak theorem there exists a
family of endomorphisms $\left\{ \widetilde{\lambda }_{i}:L\rightarrow
L_{i}\mid i\in I\right\} $, such that $\bigcap\limits_{i\in I}\ker 
\widetilde{\lambda }_{i}=\left\{ 0\right\} $. Denote $\iota
_{i}:L_{i}\hookrightarrow L$ the embedding and $\lambda _{i}=\iota _{i}%
\widetilde{\lambda }_{i}$. We have $\lambda _{i}:L\rightarrow L$ and\ $%
\bigcap\limits_{i\in I}\ker \lambda _{i}=\left\{ 0\right\} $. $\mathrm{im}%
\lambda _{i}\leq L_{i}$, $L_{i}$ is a nontrivial subgroup of $L$, so $\dim 
\mathrm{im}\lambda _{i}<\dim L$ for every $i\in I$. By reason of dimension, $%
\ker \lambda _{i}\neq \left\{ 0\right\} $ for every $i\in I$.

Let $\lambda _{i}=\varphi _{i}\oplus \psi _{i}$, where $\varphi
_{i}:V\rightarrow V$, $\psi _{i}:W\rightarrow W$. If $\ker \psi _{i}=\left\{
0\right\} $ then $\ker \lambda _{i}=\left\{ 0\right\} $ (the intersection of
a nontrivial normal subgroup of the nilpotent group with the center of group
is nontrivial - \cite[16.2.5]{KM}, similar theorem for nilpotent class Lie
algebras can be easy proved). If $\ker \varphi _{i}=\left\{ 0\right\} $ then 
$\dim \psi _{i}\left( W\right) =\dim \psi _{i}\omega \left( V,V\right) =\dim
\omega \left( \varphi _{i}\left( V\right) ,\varphi _{i}\left( V\right)
\right) =\dim \omega \left( V,V\right) =\dim W$. Hence $\ker \psi
_{i}=\left\{ 0\right\} $. It is a contradiction. So, if $\ker \lambda
_{i}\neq \left\{ 0\right\} $, then $\ker \psi _{i}\neq \left\{ 0\right\} $
and $\ker \varphi _{i}\neq \left\{ 0\right\} $. If $x\in \ker \varphi _{i}$,
then for every $v\in V$ we have $\psi _{i}\omega \left( x,v\right) =\omega
\left( \varphi _{i}\left( x\right) ,\varphi _{i}\left( v\right) \right) =0$.
So $x\in \ker \psi _{i}\omega $, i.e., skew symmetric bilinear mapping $\psi
_{i}\omega $\ is singular. $\bigcap\limits_{i\in I}\ker \psi _{i}\subset $ $%
\bigcap\limits_{i\in I}\ker \lambda _{i}=\left\{ 0\right\} $.
\end{proof}

Now for every algebra $L$ we will construct a specific geometrically
indecomposable algebra $E\left( L\right) $, such that $L\subset E\left(
L\right) $. Let $L=V\oplus W$, $\left\{ v_{1},\ldots ,v_{n}\right\} $ be a
basis of $V$, $\left\{ w_{1},\ldots ,w_{m}\right\} $ be a basis of $W$.
Then, we construct the $E\left( L\right) $ this way: $Z\left( E\left(
L\right) \right) =W\oplus T$, where $T=\mathrm{Sp}\left\{ t\right\} $ is $1$%
-dimensional $\mathbf{%
%TCIMACRO{\U{211a} }%
%BeginExpansion
\mathbb{Q}
%EndExpansion
}$-linear spaces, $E\left( L\right) /Z\left( E\left( L\right) \right) \cong
U\oplus V$, where $U$ is a $n$-dimensional $\mathbf{%
%TCIMACRO{\U{211a} }%
%BeginExpansion
\mathbb{Q}
%EndExpansion
}$-linear spaces with the basis $\left\{ u_{1},\ldots ,u_{n}\right\} $. If
the skew symmetric bilinear mapping $\omega _{L}$ is defined by $m$ skew
symmetric matrices of the size $n\times n$: 
\begin{equation*}
A^{\left( 1\right) },\ldots ,A^{\left( m\right) }
\end{equation*}%
then the skew symmetric bilinear mapping $\omega _{E\left( L\right) }$
define by $m+1$ skew symmetric matrices of the size $2n\times 2n$:%
\begin{equation*}
\left( 
\begin{array}{cc}
0 & I_{n} \\ 
-I_{n} & 0%
\end{array}%
\right) ,\left( 
\begin{array}{cc}
0 & 0 \\ 
0 & A^{\left( 1\right) }%
\end{array}%
\right) ,\ldots ,\left( 
\begin{array}{cc}
0 & 0 \\ 
0 & A^{\left( m\right) }%
\end{array}%
\right)
\end{equation*}%
i.e. 
\begin{equation}
\left[ u_{i},u_{j}\right] =0,\left[ u_{i},v_{j}\right] =-\left[ v_{j},u_{i}%
\right] =\delta _{ij}t\text{ (}1\leq i,j\leq n\text{).}  \label{basis_calc}
\end{equation}

\begin{proposition}
\label{extension1}For every algebra $L$ algebra $E\left( L\right) $ is 
\textit{geometrically indecomposable.}
\end{proposition}

\begin{proof}
Let $\omega =\omega _{E\left( L\right) }:\left( U\oplus V\right) \times
\left( U\oplus V\right) \rightarrow W\oplus T$. We assume that $\psi
:W\oplus T\rightarrow W\oplus T$ is a linear mapping such that $\psi \left(
t\right) =z\neq 0$ and $x\in \ker \psi \omega $. Denote $x=\sum%
\limits_{i=1}^{n}x_{i}u_{i}+\sum\limits_{i=1}^{n}x_{n+i}v_{i}$. $\psi \omega
\left( x,u_{j}\right) =\psi \left( -x_{n+j}t\right) =-x_{n+j}z=0$, so $%
x_{n+j}=0$ for every $j\in \left\{ 1,\ldots ,n\right\} $. Hence $%
x=\sum\limits_{i=1}^{n}x_{i}u_{i}$ and $\psi \omega \left( x,v_{j}\right)
=x_{j}z=0$, so $x_{j}=0$ for every $j\in \left\{ 1,\ldots ,n\right\} $.
Therefore $x=0$ and $\ker \psi \omega =0$. So for every linear mapping $\psi
:W\oplus T\rightarrow W\oplus T$, for which we have $\ker \psi \omega \neq 0$%
, we also have $\psi \left( t\right) =0$. By Proposition \ref{decompos_cond} 
$E\left( L\right) $ is geometrically indecomposable.
\end{proof}

$U\oplus T=H$ is an ideal of the algebra $E\left( L\right) $, $\dim \left(
H\cap Z\left( E\left( L\right) \right) \right) =1$, $\dim H/\left( H\cap
Z\left( E\left( L\right) \right) \right) =n$, $E\left( L\right) /H\cong
V\oplus W\cong L$.

\begin{theorem}
\label{main}Let $L_{1}=V_{1}\oplus W_{1}$\textit{\ and }$L_{2}=V_{2}\oplus
W_{2}$ are \textit{algebras, }$\dim V_{1}=\dim V_{2}=n$, $\dim W_{1}=\dim
W_{2}=m$,\textit{. Then }$E\left( L_{1}\right) \cong E\left( L_{2}\right) $%
\textit{, if and only if }$L_{1}\cong L_{2}$.
\end{theorem}

\begin{proof}
We denote $E\left( L_{i}\right) =H_{i}\oplus L_{i}$, $H_{i}=U_{i}\oplus
T_{i} $, $T_{i}=\mathrm{Sp}\left\{ t^{\left( i\right) }\right\} $, $\left\{
v_{1}^{\left( i\right) },\ldots ,v_{n}^{\left( i\right) }\right\} $ - basis
of $V_{i}$, $\left\{ u_{1}^{\left( i\right) },\ldots ,u_{n}^{\left( i\right)
}\right\} $ - basis of $U_{i}$ ($i=1,2$).

We assume that there is an isomorphism of algebras Lie $\alpha :E\left(
L_{1}\right) \rightarrow E\left( L_{2}\right) $. $\alpha \left( H_{1}\right) 
$ is an ideal of the algebra $E\left( L_{2}\right) $. $\dim \left( H_{1}\cap
Z\left( E\left( L_{1}\right) \right) \right) =1$, so\linebreak\ $\dim \left(
\alpha \left( H_{1}\right) \cap Z\left( E\left( L_{2}\right) \right) \right)
=1$; $H_{1}\nsubseteq Z\left( E\left( L_{1}\right) \right) $ so $\alpha
\left( H_{1}\right) \nsubseteq Z\left( E\left( L_{2}\right) \right) $. First
of all, we shall prove that $\alpha \left( H_{1}\right) \subset U_{2}\oplus
Z\left( E\left( L_{2}\right) \right) $. Let $l=u+v+z\in \alpha \left(
H_{1}\right) $ ($u\in U_{2}$, $v\in V_{2}$, $z\in Z\left( E\left(
L_{2}\right) \right) $). If $v\neq 0$, then $v=\sum%
\limits_{i=1}^{n}b_{i}v_{i}^{\left( 2\right) }$, where $b_{1},\ldots
,b_{n}\in 
%TCIMACRO{\U{211a} }%
%BeginExpansion
\mathbb{Q}
%EndExpansion
$, and exists $j\in \left\{ 1,\ldots ,n\right\} $ such that $b_{j}\neq 0$.
Then $\left[ l,u_{j}^{\left( 2\right) }\right] =-b_{j}t^{\left( 2\right) }$
by (\ref{basis_calc}) and $T_{2}\subset \alpha \left( H_{1}\right) $. Also
there exists $v_{0}\in V$, such that $\left[ v,v_{0}\right] \in
W\smallsetminus \left\{ 0\right\} $, because the skew symmetric bilinear
mapping $\omega _{L_{2}}$ is a non singular. Therefore $\left[ l,v_{0}\right]
=\left[ u,v_{0}\right] +\left[ v,v_{0}\right] \notin T_{2}$ \ and $\dim
\left( \alpha \left( H_{1}\right) \cap Z\left( E\left( L\right) \right)
\right) >1$. By this contradiction we have that $v=0$ and $\alpha \left(
H_{1}\right) \subset U_{2}\oplus Z\left( E\left( L_{2}\right) \right) $.

$\alpha \left( H_{1}\right) \nsubseteq Z\left( E\left( L_{2}\right) \right) $
so there exists $l=u+z\in \alpha \left( H_{1}\right) $ ($u\in
U_{2}\smallsetminus \left\{ 0\right\} $, $z\in Z\left( E\left( L_{2}\right)
\right) $). Because $u\neq 0$, we have as above that there is $j\in \left\{
1,\ldots ,n\right\} $ such that $\left[ l,v_{j}^{\left( 2\right) }\right]
\in T_{2}\smallsetminus \left\{ 0\right\} $. But $\dim \left( \alpha \left(
H_{1}\right) \cap Z\left( E\left( L_{2}\right) \right) \right) =1$, so $%
\alpha \left( H_{1}\right) \cap Z\left( E\left( L_{2}\right) \right) =T_{2}$.

It is clear that $\left( U_{2}\oplus Z\left( E\left( L_{2}\right) \right)
\right) \cap L_{2}=\left( U_{2}\oplus W_{2}\oplus T_{2}\right) \cap \left(
V_{2}\oplus W_{2}\right) =W_{2}$, so $\alpha \left( H_{1}\right) \cap
L_{2}\subset W_{2}\subset Z\left( E\left( L_{2}\right) \right) $ and $\alpha
\left( H_{1}\right) \cap L_{2}\subset W_{2}\cap Z\left( E\left( L_{2}\right)
\right) \cap \alpha \left( H_{1}\right) =W_{2}\cap T_{2}=\left\{ 0\right\} $%
. By arguments of dimensions we have $E\left( L_{2}\right) =\alpha \left(
H_{1}\right) \oplus L_{2}$ and, because $\alpha \left( H_{1}\right) $ is an
ideal the linear mapping $E\left( L_{2}\right) /\alpha \left( H_{1}\right)
\rightarrow L_{2}$ is an isomorphism of algebras. So we have that $%
L_{1}\cong E\left( L_{1}\right) /H_{1}\cong \alpha \left( E\left(
L_{1}\right) \right) /\alpha \left( H_{1}\right) \cong E\left( L_{2}\right)
/\alpha \left( H_{1}\right) \cong L_{2}$.

Let $\lambda =\varphi \oplus \psi :L_{1}=V_{1}\oplus W_{1}\rightarrow
L_{2}=V_{2}\oplus W_{2}$ is an isomorphism of algebras. The linear mapping $%
\varphi $ is a bijection, so in the referred above bases of $V_{1}$ and $%
V_{2}$ it is presented by the invertible matrix $F=\left( f_{ij}\right)
_{i,j=1}^{n}\in GL_{n}\left( \mathbf{%
%TCIMACRO{\U{211a} }%
%BeginExpansion
\mathbb{Q}
%EndExpansion
}\right) $. We take the matrix $G=\left( g_{ij}\right) _{i,j=1}^{n}=\left(
F^{-1}\right) ^{t}$ and by this matrix and by referred above bases of $U_{1}$
and $U_{2}$ define the linear mapping $\gamma :U_{1}\rightarrow U_{2}$.
Also, we define the linear mapping $\tau :T_{1}\ni t^{\left( 1\right)
}\rightarrow t^{\left( 2\right) }\in T_{2}$. We will prove that the linear
mapping $\gamma \oplus \varphi \oplus \tau \oplus \psi :E\left( L_{1}\right)
=U_{1}\oplus V_{1}\oplus T_{1}\oplus W_{1}\rightarrow U_{2}\oplus
V_{2}\oplus T_{2}\oplus W_{2}=E\left( L_{2}\right) $ is an isomorphism of
algebras. This mapping is a bijection, so it is necessary to prove that for
every $u^{\prime },u^{\prime \prime }\in U_{1}$ and every $v^{\prime
},v^{\prime \prime }\in V_{1}$ that the $\left( \tau \oplus \psi \right) %
\left[ u^{\prime }+v^{\prime },u^{\prime \prime }+v^{\prime \prime }\right] =%
\left[ \left( \gamma \oplus \varphi \right) \left( u^{\prime }+v^{\prime
}\right) ,\left( \gamma \oplus \varphi \right) \left( u^{\prime \prime
}+v^{\prime \prime }\right) \right] $ fulfills. But it is enough to prove
this, for basis elements of $U_{1}\oplus V_{1}$. For $1\leq i,j\leq n$ we
have $\left[ \varphi \left( v_{i}^{\left( 1\right) }\right) ,\varphi \left(
v_{j}^{\left( 1\right) }\right) \right] =\psi \left[ v_{i}^{\left( 1\right)
},v_{j}^{\left( 1\right) }\right] $, because $\varphi \oplus \psi $ is an
isomorphism of algebras, $\left[ \gamma \left( u_{i}^{\left( 1\right)
}\right) ,\varphi \left( v_{j}^{\left( 1\right) }\right) \right] =\left[
\sum\limits_{k=1}^{n}g_{ki}u_{k}^{\left( 2\right)
},\sum\limits_{s=1}^{n}f_{sj}v_{s}^{\left( 2\right) }\right]
=\sum\limits_{s=1}^{n}\sum\limits_{k=1}^{n}g_{ki}f_{sj}\delta _{ks}t^{\left(
2\right) }=\delta _{ij}t^{\left( 2\right) }=\tau \left[ u_{i}^{\left(
1\right) },v_{j}^{\left( 1\right) }\right] $ by (\ref{basis_calc}), and $%
\left[ \gamma \left( u_{i}^{\left( 1\right) }\right) ,\gamma \left(
u_{j}^{\left( 1\right) }\right) \right] =0=\tau \left[ u_{i}^{\left(
1\right) },u_{j}^{\left( 1\right) }\right] $, because $\gamma \left(
u_{i}^{\left( 1\right) }\right) ,\gamma \left( u_{j}^{\left( 1\right)
}\right) \in U_{2}$. So $E\left( L_{1}\right) \cong E\left( L_{2}\right) $.
\end{proof}

So if one can resolve the problem of the classification of the nilpotent
class $2$\ finite dimension nilpotent class $2$ Lie $\mathbf{%
%TCIMACRO{\U{211a} }%
%BeginExpansion
\mathbb{Q}
%EndExpansion
}$-algebras up to geometric equivalence, then (by Proposition \ref%
{equivalent_isom} and Proposition \ref{extension1}) he can classify up to
isomorphism the algebras $E\left( L\right) $, where $L$\ is an arbitrary Lie
nilpotent class $2$ finite dimension $\mathbf{%
%TCIMACRO{\U{211a} }%
%BeginExpansion
\mathbb{Q}
%EndExpansion
}$-algebra ($\dim Z\left( \left( E\left( L\right) \right) \right) =\dim
Z\left( L\right) +1$) and by Theorem \ref{main} he can classified up to
isomorphism all nilpotent class $2$ Lie finite dimension $\mathbf{%
%TCIMACRO{\U{211a} }%
%BeginExpansion
\mathbb{Q}
%EndExpansion
}$-algebras and all nilpotent class $2$ finite rank torsion free complete
groups.

\section{Acknowledgements.}

This research was motivated by Prof. B. Plotkin. I would like to express my
gratitude to him and to Prof. S. Margolis for their constant attention to
this work. Conversations with Prof. E. Rips, Prof. Z. Sela, Prof. E.
Hrushovski, Prof. A. Mann and Dr. E. Plotkin were very useful. After
discussions with Prof. D. Kazhdan and Prof. Yu. Drozd, I paid my attention
to the researches of Prof. V. Sergeichuk and his collaborators (\cite{Serg1}%
, \cite{Serg2}, \cite{BLS}). The debates about this problem with Dr. R.
Lipyanski led to the major break in this research, and I would like to
express my sincere gratitude. I appreciate all the authors of the paper \cite%
{BLS}, which was very contributory to this research.


\begin{thebibliography}{99}
\bibitem{Ba} Ju. Bahturin, \textit{Identities in Lie algebras,} (Russian),
Nauka, Moscow, 1985.

\bibitem{BLS} G. Belitskii, R. Lipyanski, V. Sergeichuk, \textit{Problems of
classifying associative or Lie algebras and triples of symmetric or
skew-symmetric matrices are wild.}

\bibitem{Be} A. Berzins \textit{Geometrical equivalence of algebras,}
International Journal of Algebra and Computations, \textbf{11:4} (2001), pp.
447 -- 456.

\bibitem{GrSe} F. Grunewald, D. Segal, \textit{Reflections on the
classification of torsion free nilpotent groups,} Group theory, essays for
Philip Hall, edited by K.W. Gruenberg and J.E. Roseblade, Academic Press,
1984, pp. 121 - 158.

\bibitem{GrSeSt} F. Grunewald, D. Segal, L. Sterling, \textit{Nilpotent
groups of Hirsh length six,} Mathematische Zeitschrift, \textbf{179},
(1982), pp. 219-235.

\bibitem{GrSch} F. Grunewald, R. Scharlau, \textit{A note of\ finite
generated torsion-free nilpotent groups of class 2,} Journal of Algebra, 
\textbf{58} (1979), pp. 162-175

\bibitem{KM} M. Kargapolov, Ju. Merzljakov \textit{Fundamentals of the
Theory of Groups,} Springer Verlag, New York, 1979.

\bibitem{Pl1} B. Plotkin \textit{Algebraic logic, varieties of algebras and
algebraic varieties,} Proc. Int. Alg. Conf., St. Petersburg, 1995, Walter
Gruyter, New York, London, (1999), pp. 189 -- 271.

\bibitem{Pl2} B. Plotkin \textit{Varieties of algebras and algebraic
varieties. Categories of algebraic varieties,} Siberian Advances in
Mathematics, \textbf{7(2)}, (1997), pp.64 -- 97.

\bibitem{Sch} R. Scharlau, \textit{Paare alternierender Formen,}
Mathematische Zeitschrift, \textbf{147}, (1976), pp. 13-19.

\bibitem{Serg1} V. Sergeichuk, \textit{Classification of metabelian p-groups,%
} in: Matrix problems, Inst. Mat. Ukrain. Akad. Nauk, Kiev, 1977, pp.
150-161 (in Russian); MR 58 \#11109.

\bibitem{Serg2} V. Sergeichuk, \textit{Classification problems for systems
of forms and linear mappings,} Math. USSR, Izvestiya, \textbf{31(3)},
(1988), pp. 481--501.

\bibitem{Ts} A. Tsurkov. \textit{Geometrical equivalence of nilpotent
torsion free groups.} http://arxiv.org/abs/math.GR/0411313
\end{thebibliography}
\end{document}